\documentclass{amsart}
\usepackage{srcltx}
\usepackage{mathrsfs}
\usepackage{bbm}

\usepackage{hyperref}

\def\by{\rm}
\def\paper{. \rm}
\def\book{. \rm}
\def\jour{. \it}
\def\vol{, \bf}
\def\yr{ \rm}
\def\publaddr{. \rm}
\def\publ{: \rm}
\def\pages{, \rm}

\usepackage{graphicx}

\theoremstyle{slanted}

\newtheorem{lemma}{Lemma}
\newtheorem{theorem}{Theorem}

\theoremstyle{definition}
\newtheorem{definition}{Definition}

\newtheorem{hypot}{Hypothesis} % [theorem]

\def\eps{\varepsilon}
\def\la{\lambda}

\def\G{\Gamma}
\def\a{\alpha}
\def\E{\mathop{\mathsf E}}

\def\BAlg{{\mathscr B}}
\def\fF{{\mathscr F}}
\def\fP{{\mathscr P}}
\def\scpr<#1,#2>{\langle #1, #2\rangle}
\def\scPr<#1,#2>{\left\langle #1, #2\right\rangle}
\def\Set#1{{\mathbb{#1}}}
\def\df{{\boldsymbol \delta}}

\def\conv{\mathbin{*}}

\def\Maps{\colon\:}
\def\where{\colon\:}

\def\RC{R^\circ}
\def\Mult{{\mathcal M}}
\def\Ia{${\mathcal I}$}
\def\No{No.}

\begin{document}

%\date{31.07.2012}

\author{A.\,A.~Prikhod'ko}
\address{Lomonosov Moscow State University, Moscow}
\email{sasha.prihodko@gmail.com}

\title[Ergodic maps with fast correlation decay]{%
	Ergodic automorphisms with simple spectrum characterized by fast correlation decay}

\footnotetext{%
{\it Submitted to Russian Math.\ Notes (2012)}

\smallskip
The work is supported by RFFI grant \No\,11-01-00759-a and 
the grant ``Leading Russian scientific schools'' \No\,NSh-5998.2012.1. 
}
\maketitle
\bibliographystyle{amsplain}

\begin{abstract}
The existence of measure preserving invertible transformations $T$ 
on a Borel probability space $(X,\BAlg,\mu)$ with simple spectrum is established 
possessing the following rate of correlation decay 
for a dense family of functions $f \in L^2(X,\mu)$:
\begin{equation*} 
	%\label{e1} 
	\forall\,\eps > 0 \qquad \scpr<f(T^k x),f(x)> = O(|k|^{-1/2+\eps}). 
\end{equation*} 
% для любого ${\eps > 0}$. 
% It is worth to note that, 
According to identity 
${ \scpr<f(T^k x),f(x)>  = \Hat\sigma_f(k)}$, where $\sigma_f$ denotes the spectral measure associated with~$f$, 
the rate of decay of the Fourier coefficients $\Hat\sigma_f(k)$, observed for the class of transformations introduced in the paper, 
is the maximal possible for singular Borel measures on~$[0,1]$. 

\smallskip
This note\footnotetext{Submitted to Russian Math.\ Notes.} summarizes the results of the paper 
\href{http://arxiv.org/abs/1008.4301}{arXiv:1008.4301}.

%\smallskip
%The work is partially supported by 
%	RFFI grant \No\,11-01-00759-a and 
%	Leading Russian scientific school grant \No\,NSh-5998.2012.1.
\end{abstract}

Let us consider an invertible measure preserving transformation $T$ on a Borel probability space $(X,\BAlg,\mu)$ 
and recall a question which is well-known in the spectral theory of ergodic dynamical systems 
and goes back to Banach: 
{\it Does there exist a transformation with invariant probability measure 
having Lebesgue spectrum of multiplicity one?} 

Ulam~\cite{prikh:1} (ch.~VI, \S\,6) states this problem in the following way. 

{\it Does there exist a function $f \in L^2(X,\mu)$ and a measure preserving invertible 
transformation ${T \Maps X \to X}$ ({\it an automorphism\/}), 
such that the sequence of functions $\{f(T^k x) \where k \in \Set{Z}\}$ 
is a complete orthogonal system in the Hilbert space $L^2(X,\BAlg,\mu)$? } 

One can easily construct an example of dynamical system of such kind on the space with {\em infinite\/} measure. 
Let us consider the set of integers $\Set{Z}$ as a phase space~$X$ with the standard counting measure
${ \nu(\{j\}) \equiv 1 }$, and let us define $T \Maps j \mapsto j+1$. 
Then the functions $T^k \df_0$ constitute a basis in $L^2(\Set{Z},\nu)$, 
where $\df_0(j) = 1$ if ${j = 0}$ and $\df_0(j) = 0$ if ${j \not= 0}$. 

In the class of finite measure preserving transformations the hypothesis of Banach is still open. 
Kirillov \cite{prikh:2} states a~generalized hypothesis for general Abelian group actions 
on a~space with a~finite probability measure. 

%для~каких топологических абелевых групп $G$ существуют 
%непрерывные действия группы~$G$ с конечной инвариантной мерой, 
%обладающие простым или конечнократным спектром, 
%причём мера максимального спектрального типа эквивалентна 
%мере Хаара на дуальной группе~$\hat G$? 

% In the scope of this paper 
We concentrate on the case $G = \Set{Z}$. 
For a~survey of constructions and results 
in the spectral theory of ergodic dynamical systems 
the reader can refer to \cite{prikh:3} and~\cite{prikh:4}. 
The dual group $\Hat{\Set{Z}}$ is isomorphic to the unit circle $S^1$ in the complex plane 
%
% endowed with the multiplication $(z,w) \mapsto zw$. 
% which has multiplication operation and characters in the group $\Set{Z}$ have the next form $\chi_z(k) = z^k$, 
% where ${z \in \Set{C}}$, ${|z| = 1}$. 
%
and the Fourier coefficients of any Borel probability measure $\sigma$ on~$S^1$ 
are recovered from the identity 
$$
	\Hat \sigma(k) = \int_{S^1} z^k \,d\sigma, \qquad k \in \Set{Z}. 
$$

\begin{definition}
Le $\kappa(\sigma)$ denote the liminf of the values ${\alpha \in \Set{R}}$ 
satisfying the estimate ${\Hat\sigma(k) = O(|k|^{\alpha+\eps})}$ for any ${\eps > 0}$. 
\end{definition}

Since the group $S^1$ is compact then any $\sigma$ satisfying 
$\Hat\sigma(k) = O(|k|^{-1/2-c})$ for some ${c > 0}$ is absolutely continuous with respect to 
the normalized Lebesgue measure $\la$ on~$S^1$. In~particular, 
${ \sigma = p(x)\la}$, where ${ p(x) \in L^1(S^1,\la) }$. 
Thus, any measure on $[0,1]$ with ${\kappa(\sigma) < -1/2}$ is absolutely continuous.

% The Banach and Kirillov questions are connected, in particular, to the investigation of decay rate of correlation functions 
%
Given a function $f \in L^2(X,\BAlg,\mu)$ consider a sequence of auto-correlations 
$$
	R_f(k) = \scpr<f(T^k x), f(x)>, 
$$
and the spectral measure $\sigma_f$ associated with $f$ and given by 
$
	\Hat\sigma_f(k) = R_f(k). 
$
Whenever $\kappa(\sigma_{f_j}) < -1/2$ holds for a~dense family $\{f_j\} \in L^2(X) = L^2(X,\BAlg,\mu)$, 
the spectrum of $T$ is absolutely continuous. We will prove that an extreme  
value ${\kappa(\sigma_f) = -1/2}$ (on~a~dense set of functions) is~achieved 
in the class of ergodic transformations with simple spectrum.  

\begin{theorem}\label{th1}
There exists an automorphism $T$ on a Borel space $(X,\BAlg,\mu)$ 
with simple spectrum such that $\kappa(\sigma_f) \le -1/2$ for a dense set of functions ${f \in L^2(X)}$. 
\end{theorem}

Let us define 
$$
  \kappa(T) = \inf_{\text{$\fF$ dense in $L^2(X)$}} \;\;\; \sup_{\text{$f \in \fF$}} \kappa(\sigma_f). 
$$ 
% the infimum of $\kappa^*$ such that $\kappa(\sigma_{f_j}) \le \kappa^*$
Thus, theorem~\ref{th1} states the existence of $T$ with simple spectrums and $\kappa(T) \le -1/2$. 
Observe that ${\kappa(T) = -\infty}$ for any $T$ with Lebesgue spectrum. 

% There following dichotomy theorem holds in the class of transformations with ${\kappa(\sigma_f) \le -1/2}$. 

\begin{theorem}\label{dich1}
Let $T$ be an automorphism satisfying the equality $\kappa(T) = -1/2$, 
and let $\sigma$ be the maximal spectral type of~$T$. 
Then $\sigma \conv \sigma \ll \lambda$, and, furthermore, 
the spectrum of $T$ either contains an absolutely continuous component or is purely singular, 
and for any %singular 
spectral measure $\sigma_f$ we have $\kappa(\sigma_f) = -1/2$. 
\end{theorem}

Throughout this paper we call singular Borel measures satisfying $\kappa(\sigma_f) = -1/2$ {\it Salem--Schaeffer measures}. 
This class of probability distributions were studied in the works of Schaeffer, Salem, Sigmund, Ivashev-Musatov et al.\ 
(see \cite{prikh:5,prikh:6,prikh:7}). 

The main idea of this work is to show that Salem--Schaeffer measures are found among 
spectral measures of ergodic dynamical systems. 
We propose a construction of a class of automorphisms that will serve an example of such kind. 

% We will consider two equivalent constructions: the first one applies symbolic dynamic language while the second one based on geometric construction of transformation. Then we will present the proof scheme of theorem~\ref{th1}.     

%\medskip
\begin{definition}({\it Symbolic construction})\label{def:symbolic} 
Let $\Set{A}$ be a finite alphabet and let $w_0$ be a finite word in~$\Set{A}$ 
containing at least two different letters.
Denote by $\rho_\alpha(w)$ the {\it cyclic shift\/} of $w$ to the left: % (влево) 
$$
	\rho_1(au) = ua, \qquad \rho_\alpha(u) = \rho_1^\alpha(u), \qquad 
	\text{$a \in \Set{A}$ ---  a letter}, \quad \text{$u$~--- a word}. 
$$
% где $a \in \Set{A}$ --- буква, $u$~--- слово. 
Let us construct the sequence of words $w_n$ applying the next rule: 
\begin{equation}
	w_{n+1} = \rho_{\alpha_{n,0}}(w_n)\, \rho_{\alpha_{n,1}}(w_n)\, \dots\, \rho_{\alpha_{n,q_n-1}}(w_n). 
	\label{r1}
\end{equation}
Here the sequences $q_n \in \Set{N}$ and $\alpha_{n,j} \in \Set{Z}/h_n\Set{Z}$, $h_n = |w_n|$ 
serve as parameters of the construction, and $|w|$ denotes the length of the word~$w$. 
Without loss of generality, one can assume the first entry of $w_n$ inside the bigger word $w_{n+1}$ is not touched, 
${\alpha_{n,0} \equiv 0}$. 
Then every $w_n$ is a prefix of the successor word~$w_{n+1}$, and we can define a unique infinite word $w_\infty$ 
expanding every word~$w_n$ to the right. 

Further, applying a standard procedure let us define the minimal compact subset $K \subset \Set{A}^{\Set{Z}}$  
containing all the shifts of the word $w_\infty$.    
The left shift transformation $T \Maps (x_j) \mapsto (x_{j+1})$ 
provides a topological dynamical system acting on the set $K$.
Futher, let us endow the set $K$ with a natural Borel measure $\mu$ invariant under~$T$. 
We~define the probability $\mu(u)$ of the word $u$ to be the asymptotic frequency of observing $u$ as subword in~$w_\infty$. 
% The measure $\mu$ is preserved by the transformation~$T$.     
\end{definition}

The construction and ergodic properties of the ergodic system $(T,K,\BAlg,\mu)$ are discussed in details in~\cite{prikh:8}.    
The complexity characteristics of the topological system $(T,K)$, %(без структуры измеримого пространства), 
as well as the infinite word $w_\infty$ %presenting a basis for definition are very interesting question in dynamical systems complexity problems 
are studied in ~\cite{prikh:9}. 
An infinite sequence of concateneted cyclic shifts of a fixed word % is known as {\it waltz of infinite order} 
was previously studied in the theory of recursive functions~\cite{prikh:10}. 
Setting $h_1 = 2$, $q_n \equiv 2$, $\rho_{n,0} = 0$, $\rho_{n,1} = h_n/2$, 
we see that the classical {\it Morse automorphism} (see \cite{prikh:11}) is includeded in the class definied above. 
It can be also shown that the constructed systems possess adic representation. 

%\medskip
\begin{definition}({\it Algebraic construction})\label{def:algebraic} 
Let $h_n$ be a sequence of 
positive integers such that ${h_{n+1} = q_n h_n}$, $q_n \in \Set{N}$, ${q_n \not= 1}$. 
Consider then a sequence of 
embedded lattices ${\G_n = h_n\Set{Z}}$, where  ${\G_{n+1} \subset \G_n}$,  
and the corresponding %sequence of 
homogeneous spaces ${M_n = \Set{Z}/\G_n = \Set{Z}_{h_n}}$. 
% , endowed with the natural projections ${\pi_n \Maps y + \G_{n+1} \mapsto y + \G_n}$. 

Let us fix %arbitrary 
projections $\phi_n \Maps M_{n+1} \to M_n$ definied by 
$$
	\phi_n \Maps jh_n + k \mapsto k + \alpha_{n,j} \pmod{h_n}, \qquad 
	0 \le k < h_n, \quad j = 0,1,\dots,q_n. 
$$ 
Evidently, %the projections 
$\phi_n$ preserve normalized Haar measures $\mu_n$ on the Borel spaces $M_n$. 
Define the phase space $X$ as inverse limit of spaces $(M_n,\BAlg_n,\mu_n)$, namely, set 
$$
	X = \bigl\{ x = (x_1,x_2,\dots,x_n,\dots) \where \phi_n(x_{n+1}) = x_n \bigr\}. 
$$
The measures $\mu_n$ become Borel measures $\mu$ on $X$. 
% 
%Finally, 
Let us define the transformation $T$ on the space $(X,\BAlg,\mu)$ as follows.
Any projection $\phi_n$ almost commutes with the shift transformation on $M_n$, 
$$
	\mu \bigl\{ x \where \phi_n(S_{n+1}(x_{n+1})) \not= S_n(\phi_n(x_{n+1})) \bigr\} \le h_n^{-1}, 
$$
hence, applying Borrel--Cantelli lemma, we see that $\mu$-almost surely 
the equality 
${ \phi_n(S_{n+1}(x_{n+1})) = S_n(\phi_n(x_{n+1})) }$ holds for ${n \ge n^*(x)}$, where $n^*(x)$ is a measurable function. 
Set % $T$ on the point~$x$ by the next equations:
$$
	(Tx)_n = x_n + 1, \quad n \ge n^*(x), \quad \text{and} \quad (Tx)_n = \phi_n(Tx_{n+1}), \quad n < n^*(x). 
$$
% Используя это наблюдение и применяя лемму Бореля--Кантелли, легко проверить, что 
% для $\mu$-почти всех точек ${x \in X}$ указанное соотношение коммутирования выполнено 
% The next lemma is true (see~\cite{prikh:8}). 
\end{definition}

\begin{lemma}[see~\cite{prikh:8}]
The map $T$ is a measure preserving invertible transformation %(automorphism) 
of the probability space $(X,\BAlg,\mu)$. 
\end{lemma}

The equivalence of the two constructions introduced above is verified via 
coding $T$-orbits by words $w_n$ induced by functions ${M_{n_0} \to \Set{A}}$ for some~$n_0$. 

\medskip

%\begin{constr}
In order to prove theorem~\ref{th1} we consider a certain stochastic family 
of dynamical systems $(T,X,\BAlg,\mu)$ constructed above, depending on random parameters. 
Then we show that $T$ has simple spectrum and the required rate of correlation decay 
almost surely with respect to the probability on the set of parameters.
%\end{constr}

\begin{theorem}\label{th2}
There exists a sequence $q_n \in \Set{N}$ such that the transformation $T$ 
defined above with $\alpha_{n,j}$ independent and uniformly distributed on $M_n$ 
has simple spectrum and satisfy the inequality $\kappa(T) \le -1/2$. 
% the decay of correlations 
% $\scpr<f(T^t x),f(x)> = O(|t|^{-1/2+\eps})$ for all ${\eps > 0}$ 
% for a dense family of elements ${f \in L^2(X,\mu)}$. 
\end{theorem}

%распределённых независимо и равномерно на $M_n$, автоморфизм $T$, построенный в соответствии 
%с рассмотренной конструкцией, почти наверное 
%имеет простой спектр и скорость убывания корреляций 
%$\scpr<f(T^t x),f(x)> = O(|t|^{-1/2+\eps})$ при всех ${\eps > 0}$ 
%для плотного семейства элементов ${f \in L^2(X,\mu)}$. 
%\end{theorem}

\medskip
{\it Proof.} %  of theorem~\ref{th2} %% (scheme)
The detailed proof of the simplicity of spectrum is given in~\cite{prikh:8}. 
% So, finishing the proof we would like to mention that 
% Let us mention that 
% The key ingredient in this proof is % 
It is based on the following lemma. 

\begin{lemma}[see \cite{prikh:12}]
\label{lemKS} 
Let $U$ be a unitary operator in a separable Hilbert space $H$, 
let $\sigma$ be the measure of maximal spectral type 
and let $\Mult(z)$ denote the multiplicity function of the operator~$U$. 
If ${\Mult(z) \ge m}$ on a set of positive $\sigma$-measure 
then there exist $m$ orthogonal elements of unit length $f_1,\dots,f_m$ 
such that for any cyclic space $Z \subset H$ (with respect to~$U$) 
and for any $m$ elements $g_1,\dots,g_m \in Z$ of the same length ${\|g_i\| = a}$ 
the inequality inequality holds 
$$
	\sum_{i=1}^m \| f_i - g_i \|^2 \ge m(1 + a^2 - 2a/\sqrt{m}). 
$$
\end{lemma}

%
% We start estimate decay rate of correlations. 
%
% Now we pass to the second part of the proof. 
%
In order to prove the second statement of the theorem 
it is enough to estimate the decay of correlations for $\BAlg_{n_0}$-measurable (cylinfric) functions 
$f(x)$ which are dense in~$L^2(X)$. 
Any such function $f(x)$ can be represented in the form ${f(x) = f_{n_0}(x_{n_0})}$, 
where ${ f_{n_0} \Maps M_{n_0} \to \Set{C} }$ and $x_n$ is the $n$-th coordinate of a~point~$x$. 
% (see~\cite{prikh:8}). 
Then for any $n > n_0$ 
$$
	f(x) = f_n(x_n), \quad \text{where} \quad f_{n+1}(x_{n+1}) = f_n(\phi_n(x_{n+1})). 
$$
Given a function $f(z)$ with zero mean define the {\it cyclic correlations\/} 
$$
	\RC_n(t) = \int_{M_n} f_n(j+t)\,\overline{f_n(j)} \,d\mu_n(j), 
$$
where $j+t$ is the sum in the group $M_n$, i.e.\ $j+t \pmod{h_n}$. 
Taking into account the conditions on the distribution of the randomparameters $\alpha_{n,j}$ 
it can be easily shown that $\mu$-a.s.\ ${\RC_n(t) \to R_f(t)}$ for any~$t$, hence, the distributions 
$\Hat\RC_n$ converges weakly to the spectral measure $\sigma_f$ (but we need only the first convergence). 
Now let us consider (the most important) case, when ${t = sh_n}$, ${s \not= 0}$. 
For this special value of the argument $t$ the following recurrent identity holds
$$
	\RC_{n+1}(t) = \frac1{q_n} \sum_{k=0}^{a_n-1} \RC_n(\a_{n,k+s}-\a_{n,k}). 
$$
Applying expectation operator with respect to the probability on the parameters' space, we obtain
$\E \RC_{n+1}(t) = \E \RC_n(\a_{n,k+s}-\a_{n,k}) = 0$ and 
$$
	\E |\RC_{n+1}(t)|^2 = \frac1{q_n^2} 
		\E \sum_{k,\ell=0}^{q_n-1} \RC_n(\a_{n,k+s}-\a_{n,k})\,\overline{\RC_n(\a_{n,\ell}-\a_{n,\ell+s})}. 
$$
Observe that all the terms in the above sum are zero except the terms such that $\{k,k+s\} = \{\ell,\ell+s\}$. 
If $h_n$ is odd these are always the terms with ${k = \ell}$, and if $n$ is even we should also count 
the terms given by ${k+s = \ell}$, and ${\ell+s = k}$. Clearly, the latter contributes $O(q_n^{-1})$ to the sum 
over all~$s$, so without loss of generality we can assume that all $h_n$ are odd. We have then  
$$
	\E |\RC_{n+1}(t)|^2 = \frac1{q_n} \E |\RC_n(\a_{n,k+s}-\a_{n,k})|^2, 
$$
hence, $\E |\RC_{n+1}(t)|^2 = h_{n+1}^{-1} \E\|\RC_n\|^2$, 
where $\|\cdot\|$ is a standard form in $L^2(\Set{Z})$, and the function $\RC_n$ is restricted to $[0,h_n-1]$. 
The same arguments slightly modified lead the equality $\E |\RC_{n+1}(t)|^2 = h_{n+1}^{-1} \E\|\RC_n\|^2$ 
for any $t \in (h_n,h_{n+1})$. 
Thus, one can see that 
$$
	\E \|\RC_{n+1}\|^2 \le 2 \cdot \E\|\RC_n\|^2. 
$$
It follows that % further, applying a set of additional estimations we obtain
$\RC_{n+1}(t) = O(|t|^{-1/2+\eps})$ for any $\eps > 0$ 
the second statement of the theorem is verified.~$\Box$

\begin{definition}
We say that an automorphism $T$ of a Borel probability space $(X,\BAlg,\mu)$ 
admits {\it approximation of type \Ia\/} if for any finite partition $\fP$ and any $\eps > 0$ 
there exists a subset $\Omega_\eps \subset X$ of the measure $1-\eps$ and a word $W_\eps$ 
such that for all $x \in \Omega_\eps$ the infinite word generated by $\fP$-coding of the $x$-orbit 
is $\eps$-covered by a sequence of words $\Tilde W_j$ which are 
$\Bar d$-$\eps$-close to cyclic shifts $\rho_{\a_j}(W_\eps)$ of the word $W_\eps$. 
\end{definition}

This property is a metric invariant. In particular, the class of maps satisfying type \Ia\ approximation includes rank one transformations. 
Clearly, the main construction of this paper generates transformations of type~\Ia. 

\begin{hypot}
Let $T$ be an automorphism constructed according to definition~\ref{def:symbolic} (or definition~\ref{def:algebraic}).  
Consider an arbitrary $\BAlg_n$-measurable function $f$ with zero mean. Then ${\kappa(\sigma_f) \ge -1/2}$. 
\end{hypot}

\begin{hypot}
Assume that an automorphism $T$ admits approximation of type~\Ia, 
and suppose that $\xi_n$ are finite partitions generating, for any fixed finite partition $\fP$, 
approximating sequence of words $W_{\eps_n}$ with ${\eps_n \to 0}$. 
% Let $\xi_n$ be the partitions 
Then 
$$  
  \liminf_{k \to \infty} \;\;\; \inf_{\text{$f$ is $\xi_n$-measurable}} \kappa(\sigma_f) \ge -1/2
$$ 
\end{hypot}

\medskip
The author thanks A.\,M.~Vershik, V.\,V.~Ryzhikov, J.-P.~Thouvenot, 
E.~Janvresse, T.~de la Rue, K.~Petersen, B.~Weiss, and the participants of seminar 
``Probability theory and statistical phyics'' for discussions and remarks.

%===============Список литературы==================

\end{document}